\input amssym.def
\input amssym.tex
\magnification 1200
\font\sectionfont=cmbx12 at 12pt

\def\min{\setminus}
\def\BZ{{ \Bbb Z}}

\def\BQ{{ \Bbb Q}}
\def\BC{{ \Bbb C}}
\def\BP{{\Bbb P}}
\def\calA{{\cal A}}
\def\calE{{\cal E}}
\def\calF{{\cal F}}
\def\calO{{\cal O}}
\def\calH{{\cal H}}
\def\na{\nabla}
\def\pr{\prime}
\def\display{\displaystyle}
\def\item{\vskip .03 in}
\def\O{\Omega}
\def\a{\alpha}

\def\g{\gamma}

\def\sec#1{\vskip 1.5pc\noindent
{\hbox 
{{\sectionfont #1}}}\vskip 1pc}

\def\theo#1#2{\vskip 1pc
  {\bf   Theorem\ #1.} {\it #2}
\vskip 1pc}
\def\lem#1#2{\vskip 1pc{\bf 
\ Lemma #1.}
  {\it #2}\vskip 1pc}
\def\cor#1#2{\vskip 1pc{\bf 
\ Corollary\ #1.}
  {\it #2}\vskip 1pc}
\def\prop#1#2{\vskip 1pc{\bf 
 Proposition \ #1.} 
{\it #2}}
\def\defi#1#2{\vskip 1pc 
{\bf Definition \ #1.}
{\it #2}\vskip 1pc }

\def\proof#1{\vskip 1pc {\bf  
Proof.} {\rm #1} 
\vskip 1pc}

\def\hook{\hookrightarrow}
\def\Sp{Springer Verlag}

\def\noblackbox{\overfullrule=0pt}
\noblackbox
\def\NSF95{\footnote*{This research was supported 
in part by NSF grant
 DMS-9504522.}}
\def\qed{\vrule height4pt 
width3pt depth2pt}

\rightline{January 30, 1997}
\centerline{\bf A correspondence dual to McKay's}
\vskip 1pc
\centerline{Jean-Luc Brylinski\NSF95}
\vskip 1pc

\sec{1. Introduction}

It is well-known from the work of DuVal [{\bf DuVal1}]
 and M. Artin [{\bf A}]that
there is a one-to-one correspondence between finite subgroups $G$
of $SU(2)$ and Coxeter-Dynkin diagrams $\Delta$ of type $A,D,E$.
This involves a minimal resolution of singularities $\tilde X$
of the singular algebraic surface $\BC^2/G$.
Around 1980 McKay found a deep correspondence between vertices
of the affine Coxeter-Dynkin diagram and irreducible representations
of the group [{\bf McKay1}] [{\bf Mckay2}]. Several  systematic 
representation-theoretic proofs were
given by Kostant [{\bf Ko}], Steinberg [{\bf St2}], Springer [{\bf Sp}].
A geometric interpretation of the correspondence
 was given by Gonzalez-Sprinberg
and Verdier [{\bf GS-V}] and also by Kn\"orrer [{\bf Kn}].

There is also a dual correspondence, this time between
vertices of $\Delta$ and non-trivial conjugacy classes of $G$.
This dual correspondence was introduced by Ito
and Reid [{\bf I-R}] in the more general context
of a finite subgroup of $SU(n)$.
The construction is in fact very simple topologically:
we interpret $G$ as the fundamental group of the 
complement of the exceptional
divisor in $\tilde X$. Then each vertex of $\Delta$ corresponds
to a component of this divisor, and there is an associated class
of a small loop encircling said component. The dual correspondence
was studied by Ito and Reid from the point of view
of valuations on function fields. From the description
of the fundamental group due to Mumford [{\bf Mu}] one sees
that this gives a one-to-one correspondence between vertices
of the diagram and non-trivial conjugacy classes. This result
amounts to the dimension $2$ case of a more general theorem
proved by Ito and Reid [{\bf I-R}]. The dual
correspondence
has a number of interesting further properties, which are detailed
in Theorem 4.1. These properties involve the three (or two)
so-called special conjugacy classes corresponding to ends
of the diagram, and
the description of the edges of the diagram involve pairs of commuting
elements $x,y$ such that $y$ is special; then the conjugacy classes
of $x$ and $xy$ are joined by an edge.

There seem to be intriguing connections between the McKay correspondence
and the dual correspondence. We prove a determinantal formula
concerning an element $g_j$ of $G$ associated to a vertex $v_j$
of $\Delta$ and the irreducible representation $E_k$ associated
to a vertex $v_k$:
$$det(g_j,E_k)=exp(-2\pi i(C^{-1})_{jk}),$$
where $C^{-1}$ is the inverse of the Cartan matrix. 
One tool we use in proving this formula is the geometric
description of the McKay correspondence
in [{\bf GS-V}] by means of the first Chern class
of the vector bundle associated to a representation
of $G$.   We also use the properties
of  vector bundles with integrable
connections admitting logarithmic poles, in particular
the computation of the first Chern class
from the residues of the connection.

The paper ends with some remarks  on the
matrix-valued Fourier transform which results from comparing the
two correspondences.

This work was written in May-June 1996 as I was
visiting Harvard University. I thank David Kazhdan and the Harvard
Mathematics Department
for their hospitality. I am grateful to H\'el\`ene Esnault for correspondence
and information about her work [{\bf E-V}] and for her
useful remarks on a first draft of this paper. I thank Victor Batyrev
for pointing out to me the paper [{\bf I-R}]. I am grateful to Igor Dolgachev
and to John McKay for reading a first version of this paper and making
many valuable comments.

\sec{2. The McKay correspondence}
There is a by now classical correspondence between
conjugacy classes of finite subgroups of $SU(2)$
(or equivalently, of $SL(2,\BC)$)
and  simply-laced Coxeter-Dynkin diagrams
(thus of type $A_n$, $D_n$ or $E_n)$. The correspondence
was constructed by DuVal [{\bf DuVal1}] using algebraic
geometry. It may be phrased as follows.
Given a finite subgroup $G\subset SU(2)$ one can construct
the singular algebraic surface $X=\BC^2/G$, quotient
of the affine plane by the action of $G$.
Let $f:\BC^2\to X$ be the quotient map;
the point $o=f(0)$ is called the origin of $X$. The singular
locus of $X$ is reduced to $o$, unless $G=1$, in which
case $X=\BC^2$ is smooth.
$X$ is a normal affine surface, whose algebra
of regular functions is the algebra $\BC[x,y]^G$
of invariants in the polynomial algebra on two generators.
There is a minimal resolution of singularities
$p:\tilde X\to X$, so 

- $\tilde X$ is a smooth algebraic
surface;

- $p$ is a proper regular mapping, which induces
an isomorphism over the open subset $U=X\min \{ o\}$ of $X$;

- the minimality of the resolution means that
 $\tilde X$ does not contain
any rational curve of self-intersection $-1$.

Then the reduced fiber $D=p^{-1}(o)_{red}$ is a curve, which
is a union of smooth rational curves $D_1,\cdots,D_r$.
Any two  curves $D_i$ and $D_j$
intersect transversally at at most one point of $\tilde X$.
One can then construct a graph $\Delta$ such  that the vertices
$v_1,\cdots,v_r$
of $\Delta$ correspond to the curves $D_j$, and where
one joins the vertices $v_i$ and $v_j$ by an edge whenever
the divisors $D_i$ and $D_j$ intersect.

One associates to the system of curves $(D_j)$ a  square matrix
$A=A(G)$ of size $r$, called the intersection matrix,
 such that $A_{ij}$ is the intersection number
$D_i\cdot D_j$ on the smooth surface $\tilde X$.
On the other hand the graph $\Delta$ has an incidence matrix $M$.
Since any $D_i$ has self-intersection $-2$ we have
$A=-2 Id+M$.

We then have

\theo{2.1}{(Du Val [{\bf DuVal1}], M. Artin [{\bf A}]) 
(1) For any finite subgroup $G$ of $SU(2)$, the graph
$\Gamma$ is  a simply-laced Coxeter-Dynkin diagram.
\item (2) The Cartan matrix $C$ of the Coxeter-Dynkin diagram
$\Delta$ is the opposite of the intersection matrix $A$.
\item (3) This construction induces a one-to-one correspondence
between  conjugacy classes of finite subgroups of $SL(2,\BC)$
and simply-laced Coxeter-Dynkin diagrams.}

We give the table showing the simply-laced diagrams and the
corresponding finite subgroups of $SU(2)$.
For a regular polyhedron, we have the corresponding
symmetry group $H\subset SO(3)$. Its inverse image
$G$ in the double cover $SU(2)$ is the corresponding
 binary polyhedral group.

\vskip 1.5pc
\centerline{{\bf Table of subgroups of $SU(2)$}}
\vskip 1pc

\vbox{\tabskip=0pt \offinterlineskip
\def\tablerule{\noalign{\hrule}}
\halign to345pt{\strut#& \vrule#
\tabskip=1em plus2em&
\hfil#& \vrule#&\hfil#\hfil& \vrule#&
#\hfil& \vrule#\tabskip=0pt \cr\tablerule
&&\multispan5\hfil FINITE SUBGROUPS OF 
$SU(2)$\hfil&\cr\tablerule
&&\omit\hidewidth $\Delta$\hidewidth&&
\omit\hidewidth order of $G$\hidewidth&&
\omit\hidewidth $G$\hidewidth&\cr\tablerule
&&$A_n$&& $n+1$&&cyclic&\cr\tablerule
&&$D_n,n\geq 3$&&$4n-8$&&binary dihedral&\cr\tablerule
&&$E_6$&&$24$&&binary tetrahedral&\cr\tablerule
&&$E_7$&&$48$&&binary octahedral&\cr\tablerule
&&$E_8$&&$120$&&binary icosahedral&\cr\tablerule}}
\vskip 2pc

The McKay correspondence on the other hand involves
simply-laced affine Coxeter-Dynkin diagrams
 [{\bf McKay1}] [{\bf McKay2}].
Given a simply-laced Coxeter-Dynkin diagram $\Delta$, there
is a corresponding affine diagram
 $\tilde\Delta$, which is obtained
by adding one vertex $v_0$ to $\Delta$. We need
 to explain for which
$i\in\{ 1,\cdots,r\}$ the vertex $v_0$ and $v_i$
 are linked by an edge.
This requires introducing the root system $R$ 
corresponding to $\Delta$.
This a finite subset of an euclidean vector
 space $E$ of dimension $r$,
consisting of vectors of length $2$.
The vertices $v_1,\cdots,v_r$ correspond to the simple roots
$\a_1,\cdots,\cdots,\a_r$, with respect to system $R_+$
of positive roots. The Cartan matrix
is given by $C_{ij}=(\a_i,\a_j)$. There is a
 longest root $\theta$ (so that $\theta$
is a positive root, and $\theta+\a_i$ is not
 a root for $j=1,\cdots,r$).
Then the new vertex $v_0$ of $\tilde\Delta$ corresponds to 
$\a_0:=-\theta$. There is an edge in $\tilde\Delta$
between the vertices $v_0$ and $v_i$ if and only
 if $(\a_0,\a_i)\neq 0$.

Each vertex $v_i$ of $\tilde\Delta$ is labeled
 by a positive integer $m_i$
in such a way that $m_0=1$ and
 $\sum_{j=10}^r~m_i\a_i=0$. Equivalently
we have $\theta=\sum_{j=1}^r~m_j\a_j$.

We can now state the McKay correspondence.

\theo{2.2}{(Mckay [{\bf McKay}]) Let $G$ be a finite subgroup
of $SL(2,\BC)$ and let $\tilde\Delta$ be
 the corresponding affine
Coxeter-Dynkin diagram.
(1) There is a one-to-one correspondence
$i\mapsto R_i$ between vertices of $\tilde\Delta$
and equivalence classes of irreducible representations
of $G$. The dimension of $R_i$ is equal to $m_i$.
\item(2) Let $E$ be the two-dimensional representation
of $G$ in $\BC^2$. Then for any $i\in \{ 0,\cdots,r\}$ we have
$$R_i\otimes E\tilde{\to}\oplus_{j~
{\rm incident~to}~i}~R_j\eqno(2)$$}

This correspondence was constructed empirically by McKay.
Coherent proofs were given in [{\bf St2}] [{\bf Sp}].
The representation-theoretic and invariant-theoretic
aspects of the correspondence were developed further in
[{\bf Ko}].
A geometric construction was given by Gonzalez-Sprinberg
and Verdier [{\bf GS-V}]; this will be used in \S 5.

\sec{3. The special conjugacy classes.}

We will use a well-known topological interpretation
of the group $G\subset SL(2,\BC)$.

\lem{3.1}{We have 
$$G\tilde{\to}\pi_1(X\min \{ o\})
=\pi_1(\tilde X\min D)\eqno(3-1)$$}

We did not specify a base point in Lemma 3.1. This is because
we only need the isomorphism
$G\tilde{\to}\pi_1(\tilde X\min D)$ up to conjugation.

\proof{The space $X\min \{ o\}$ is the
 quotient of the simply-connected
space $\BC^2\min \{ 0\}$ by the action of $G$. Because
$G\subset SL(2,\BC)$, the action of $G$ on $\BC^2\min \{ 0\}$
is fixed point free. Thus $\BC^2\min \{ 0\}\to X\min \{ o\}$
is a Galois covering with group $G$. \qed}

For any component $D_i$ of $D$ ($1\leq i\leq r$),
there is a well-defined conjugacy class  in 
$\pi_1(\tilde X\min D)$
which corresponds to a small oriented loop $\g_i$
around the divisor $D_i$. Of course the precise construction
of this loop depends on the base point but the conjugacy class
is well-defined. Using the isomorphism 3-1,
this defines a conjugacy class in $G$, 
which will be denoted by $O_i$.

\defi{3.2}{The dual McKay correspondence is the map 
$$\{ {\rm vertices~of}~\Delta\}\to G/conj$$
which maps $v_i$ to $O_i$.}

This is a more topological version of the construction
of Ito and Reid [{\bf I-R}], which is phrased
in terms of valuations and is purely algebraic
i.e., invariant under automorphisms of the field
$\BC$.

The main properties of the dual
correspondence will be given in \S 4.
These will involve some special
conjugacy classes in $G$, which are indexed by the
{\it ends} of the graph $\Delta$. There are two ends
for the graph $A_n$ and three ends
for the graph $D_n$ and $E_n$.

For this purpose we consider the induced action
of $G$ on the projective line $\BC\BP^1$ of lines
in $\BC^2$. This action factors through an effective
action of the image $H$
of $G$ in the quotient group $PGL(2,\BC)=SL(2,\BC)/\pm 1$.
The quotient $\BC\BP^1/H$ is isomorphic
to the projective line, so we have a ramified 
covering $\pi:\BC\BP^1\to \BC\BP^1/H$, which is a  Galois
covering of group $H$.
We note the well-known lemma

\lem{3.3}{(1) If $\Delta$ is of type $A_n$ for $n$
odd, we have: $G\tilde{\to}H$.
\item(2) In every other  case we have an exact sequence
$$1\to \pm \{ 1\}\to G\to H\to 1\eqno(3-2)$$}

There are three  ramification points
of $\pi$ in $\BC\BP^1/H$,
except in case $A_n$,  when there are only two.
We will give a bijection between the ramification set of $\pi$
and the set of ends of $\Delta$.

For each ramification point $q\in \BC\BP^1/H$ pick
 a point $\tilde q\in\pi^{-1}(q)$
which corresponds to a line $l\subset \BC^2$.
 We have a natural mapping
$l\hook \BC^2\to \BC^2/G$.
The inclusion
$l\min \{ 0\}\hook \BC^2\min \{ 0\}$ gives a  regular
mapping
$$l\min \{ 0\}\to X\min \{ 0\}=
\tilde X\min D.\eqno(3-3)$$ Because
the map $\tilde X\to X$ is proper, it
 follows that we can extend
the mapping (3-3) to a regular mapping
 $\phi:l\to \tilde X$.
The point $\phi(0)$ of $\tilde X$ is independent of the choice
of $\tilde q\in \pi^{-1}(q)$.
Thus to $\tilde q$ we have attached the point $x=\phi(0)$
of $\tilde X$. Let $C_q$ be the image of 
$\phi$, which only depends
on $q$, not on the choice of $\tilde q$. The map
$l\to C_q$ is a ramified Galois covering with Galois group
equal to the stabilizer $G_l$ of $l$ in $G$.

\prop{3.4}{(1) For a ramification point $q\in \BC\BP^1/H$,
the corresponding point $x$ of $\tilde X$ belongs
to only one divisor $D_j$. The curve $C_q\subset \tilde X$
is a smooth curve which meets $D_j$ transversally.
\item (2) The vertex $v_j$ is an end of the graph $\Delta$.
\item (3) The map $q\mapsto v_j$ gives a bijection between
the set $S$ of ramification points of $\pi:\BC\BP^1\to \BC\BP^1/H$
and the set of ends of $\Delta$.}

\proof{In case of a cyclic group $G$ of order $n$, the statement is
easy to prove using the explicit description of $\tilde X$
given in [{\bf Br}] or in [{\bf GS-V}]. In that case
one of the lines $l$ is the line $x=0$. There is a covering
of $\tilde X$ by $n-1$ affine open sets, each of them isomorphic
to the affine plane $\BC^2$. Consider the first  open set $U_1$
with coordinates $(u,v)$. These can be chosen so that $x=v^nu^{n-1}$
and $y=u$. Then the point $\phi(0)$ is the point $u=0,v=0$
which verifies the statement, as the divisor $v=0$ is the
divisor $D_j$ corresponding to an end of the graph.
A similar argument can be applied to the line $y=0$.
For the other cases one can use the results of Brieskorn [{\bf Br}]
to  deduce them from the cyclic case. First
 we make a preliminary observation
concerning the natural action of $\BC^*$ on $\BC^2$ by dilations,
which induces an algebraic $\BC^*$-action on $X$ and on $\tilde X$.
The punctured lines $l\min \{ 0\}\subset \BC^2$ are $\BC^*$-orbits.
Their images $C_l\min \{\phi(0)\}$ in $X$ are therefore
also $\BC^*$-orbits, and this describes all the $1$-dimensional
orbits. Among the orbits of dimension $1$,
 those corresponding to ramification
points of $\pi$ are characterized by the fact
 that the action of $\BC^*$
is not free (the $m$-th roots of unity act trivially, if $m$
is the order of $G_l$).
Now let $Y$ be the blow-up of $o\in X$. According to [{\bf Br}]
$Y$ has only isolated singularities (as many as the ends
of $\Delta$) which are rational surface singularities of
type $A_n$. The isomorphism between the germ of $Y$ at
such a point $q$ and the germ of $\BC^2/\mu_{n+1}$
at $o$ can be made $\BC^*$-equivariant. Now for the line $l$
corresponding to ramification point, the corresponding
limiting point in $Y$ is a fixed point of $\BC^*$,
so it is one of the singular. Furthermore, the corresponding
germ of orbit in the singular surface $\BC^2/\mu_{n+1}$
is a special orbit on which $\BC^*$ does not act freely.
This special orbit itself corresponding to a line $x=0$ or $y=0$.
Now the resolution of singularities $\tilde X$ is obtained
from $Y$ by minimally resolving each singular point. The effect
of this operation is already understood. \qed}

Now for any  ramification point $q\in S\subset\BC\BP^1/H$
there is a well-defined conjugacy class $V_i$ in $H$, which
 is defined as follows.  There is a group homomorphism
$f:\pi_1([\BC\BP^1/H]\min S)\to H$, which is only
well-defined up to conjugacy. Take  a small loop
in $[\BC\BP^1/H]\min S$
encircling the point $q$, and let $h_i$ be its image in $G$.
Then $V_i$ is the conjugacy class of $h_i$.
We can now state

\lem{3.5}{The conjugacy class $V_i\subset H$ is the image of the
conjugacy class $O_i\subset G$ under the canonical homomorphism
$G\to H$.}

\proof{Clearly a representative $h_i$ of $V_i$ is the image
in $H$ of  the generator of the stabilizer $G_l$ of $l$ which
admits $\display e^{2\pi i\over s}$ as  an eigenvalue, where $s$
is the order of $G_l$. This is the same as the image
 under the group homomorphism
$$\pi_1([l\min \{ 0\}]/G_l)\to \pi_1([\BC^2\min
 \{ 0\}]/G)\tilde{\to}G\to H$$
of a small loop  in $[l\min \{ 0\}]/G_l$ which 
turns once around the point
$0$.
On the other hand the curve $C_q=\phi(l)\subset
 \tilde X$ is the closure
of $[l\min \{ 0\}]/G_l\subset \tilde X\min D$.
From Proposition 3.3 
we see that the conjugacy
class $O_i\subset G$ is represented by  a 
 small loop inside this curve
which turns once around the point $x=\phi(0)$. \qed}

In case there are three ramification points
 $q_1,q_2,q_3$ we can choose
representatives $h_1,h_2,h_3$ of the three 
corresponding conjugacy classes
in $H$ such that $h_1h_2h_3=1$ (indeed this relation
holds among the conjugacy classes
in $\pi_1([\BC\BP^1/H]\min S)$
 corresponding to the three punctures).
it is natural to ask what relation exists
 among the conjugacy classes
in $G$. 

\lem{3.6}{(1) In cases $D_n$ and $E_n$,
the conjugacy classes $C_1,C_2,C_3$ corresponding
to the ends $v_1,v_2,v_3$ of the graph $\Delta$ have
representatives $g_1,g_2,g_3$ which satisfy
$$g_1g_2g_3=-1\eqno(3-4).$$
\item(2) In the case $A_n$ the conjugacy classes
$g_1$ and $g_2$ corresponding to the two ends
of $\Delta$ satisfy $$g_1g_2=1\eqno(3-5).$$ }

\proof{This is easily checked using the explicit description
of the group $G$ given for instance in [{\bf Coxeter}]
or in [{\bf DuVal2}]. \qed}

We state the following result only in the case of a graph
with three ends (the case of two ends is simpler
 and is left to the reader).

\prop{3.7}{(Coxeter) (1) The group $H$ is
 isomorphic to the abstract
group with generators $h_1$,$h_2$, $h_3$
 and defining relations
$$h_j^{m_j}=1,h_1h_2h_3=1\eqno(3-6)$$
where $m_j$ is the order of $h_j$ in $H$, which
 is also equal to the length
from $v_j$ to the central vertex $v_{cen}$.
\item(2) The group $G$ is 
isomorphic to the abstract
group with generators $g_1$,$g_2$, $g_3$ and
 defining relations
$$g_1^{m_1}=g_2^{m_2}=g_3^{m_3}=g_1g_2g_3, (g_1g_2g_3)^2=1
\eqno(3-7).$$}

Of course the central element $c=g_1g_2g_3$ is of order $2$ and
 corresponds to $-1\in SL(2,\BC)$.

\sec{4. The dual McKay correspondence.}

The dual McKay correspondence was introduced in Definition 3.2.
It associates to a vertex $v_i$ of $\Delta$ a conjugacy class
$O_i$ in $G$, which is defined topologically as the class
of a small loop in $\tilde X\min D$ encircling the divisor
$D_i$.

The main properties of this correspondence are given
in Theorem 1 below. Except in case
 $A_n$ there is a central vertex
$v_{cen}$ of $\Delta$, and there are three branches.

We will use the notion of {\it canonical
 numbering} of the vertices
of the tree $\Delta$. This means that the
 vertices are numbered
$1,\cdots,r$ and that for any $2\leq k\leq r$
 the corresponding vertices
$v_1,\cdots v_k$ are the vertices of a
 subtree. Any canonical numbering
gives an ordering of the set of vertices. 
Such an ordering will be called
canonical. This notion was used
by vonRandow [{\bf vR}].

\theo{4.1}{(1) The correspondence $v_i\mapsto C_i$ gives
a
bijection between the set of vertices of $\Delta$ and the set
of non-trivial conjugacy classes of $G$.
\item(2)  The ends of $\Delta$ correspond to the special
conjugacy classes of $G$.
\item(3) (cases $D_n$, $E_n$) The conjugacy class corresponding
to the central vertex $v_{cen}$ consists of the central element
$-1$.
\item(4) A branch $v_1,\cdots,v_m$ of $\Delta$ corresponds
to a ``geometric progression''\hfill\break
 $g,g^2,\cdots,g^m$.
\item(5) Two vertices $v_i$ and $v_j$ belong to the same branch
if and only if the corresponding conjugacy
classes $C_i$, $C_j$  have representatives $g_i$ and $g_j$ which commute.
\item (6) Two vertices $v_i$ and $v_j$ are joined by an edge
if and only there exists a representative $g_i$ 
of $C_i$ and a representative
$u$ of some special conjugacy class such that
$u$ commutes with $g_i$ and such that $ug_i$ belongs to $C_j$.
\item (7) Pick a canonical ordering of the vertices
of $\Delta$. Let $v_i$ be any vertex of $\Delta$, and let
$v_j,v_k,\cdots$ be the ordered set of neighbors of $v_i$. 
Then there are representatives $g_i,g_j,\cdots$ of 
the corresponding conjugacy classes
such that
$$g_i^2=g_jg_k\cdots.\eqno(4-1)$$
\item(8) For any canonical ordering of the
 vertices of $\Delta$, once can choose
an element $g_i$ of each conjugacy class $C_i$
such that (4-1) holds for any vertex, and such that
\item
$g_i$ and $g_j$ commute whenever the vertices $v_i$ and 
$v_j$ are joined by an edge. \hfill(4-2)\break
\item Then $G$ is described as an abstract group
as the group generated by these elements $g_i$
subject to these two types of relations.}

\proof{We have again $G=\pi_1(\tilde X\min D)$.
The fundamental group $\pi_1(\tilde X\min D)$ has been described
by Mumford [{\bf Mu}]  in terms of  precisely chosen loops
around $D_i$ with class $g_i\in G$. Mumford proved statement (8).
Now by Proposition 3.6 an end $v_j$ of $\Delta$ corresponds
to the conjugacy class of some $g_j$ whose order
is exactly the length $m_j$ of the branch which ends at $v_j$.
Index the vertices on the branch by $1,\cdots,m_j$.
Then applying (4-1) inductively, we see that the vertex labeled
by $k$, $1\leq k\leq m_j$ corresponds to the conjugacy class
of $g_j^k$. In case $A_{n-1}$, $G=\mu_n$, with the vertices
labeled linearly by $\{ 1,\cdots,n-1\}$, we see that
the vertex labeled by $k$ corresponds to
$\display e^{2\pi ik\over n}$, and all statements of the theorem
are clear. So we now assume that we are in the case $D_n$ or $E_n$,
which means there is a central vertex $v_{cen}$. Now
for $k=m_j$, we have the other end of the branch, which is the central
vertex; the corresponding conjugacy class
is represented by $g_j^{m_j}=-1$ (cf. Proposition 3.6 (2)).
At this point we can describe
 graphically the dual correspondence.
$$\matrix{&&&&&&&&y\cr
x&&x^2&&\cdots&&c&&\cr
&&&&&&&&yx}\eqno(D_n)$$
where $G$ is generated by $x$ and $y$
 with defining relations:
$$x^{n-1}=y^{n-1}=c, c^2=1, yxy^{-1}=x^{-1}\eqno(4-3)$$
\vskip 1pc
$$\matrix{x&&x^2&&c&&y^2&&y\cr
&&&&&&&&\cr
&&&&z&&&&}\eqno(E_6)$$
where $G$ is generated by $x,y,z$
with defining relations
$$x^3=y^3=z^2=c, c^2=1,  xyz=-1\eqno(4-4)$$ 
\vskip 1pc
$$\matrix{x&&x^2&&x^3&&c&&y^2&&y\cr
&&&&&&&&&&\cr
&&&&&&z&&&&}\eqno(E_7)$$
where $G$ is generated by $x,y,z$ with defining relations
$$x^4=y^3=z^2=c, c^2=1,  xyz=-1\eqno(4-5)$$
\vskip 1pc
$$\matrix{x&&x^2&&x^3&&x^4&&c&&y^2&&y\cr
&&&&&&&&&&&&&\cr
&&&&&&&&z&&&&&&}\eqno(E_8)$$
where $G$ is generated by $x,y,z$ with defining relations
$$x^5=y^3=z^2=c, c^2=1,  xyz=-1\eqno(4-6)$$
\vskip 1pc
The remaining statements can then be checked directly.
Since the number of conjugacy classes of $G$ is equal
to the number of vertices of $\Delta$, 
(1) will follow if we show that
distinct vertices $v_i,v_j$ correspond to
 distinct conjugacy classes.
This is easy to see if $v_i$ and $v_j$ are on the same branch.
Then the trace of $g^i$ in the two-dimensional representation
$\BC^2$ of $G$ is equal to $2 cos({\pi i\over m})$ and
for $i\neq j$, $1\leq i,j\leq m$ we have
$2 cos({\pi i\over m})\neq 2 cos({\pi j\over m})$.
Then there are some cases to be considered where vertices
on different branches correspond to conjugacy classes
of elements of the same order. In case 
$D_n$ it is easy to see that
$y$ and $yx$ are not conjugate. In case $E_6$ 
one checks that $x^2$ and $y^2$
are not conjugate (this is related to the fact that
 there are two distinct
conjugacy classes of rotations of order $3$ 
in the symmetry group of the
tetrahedron). This implies that $x$ and $y$
 are not conjugate.
In case $E_7$ one checks that $x^2$
 and $y$ are not conjugate, as
their images in the symmetry group of the cube are not
conjugate: the first is a flip whose axis 
 goes through the center
of two faces, the second one a flip whose axis
 goes through the middle
of two edges.
Statement (5) is easy to check, and then (6) follows directly. 
\qed}

Part (1) of the theorem is due to Ito
and Reid [{\bf I-R, Theorem 1.4}]. In fact, for an arbitary finite
subgroup $G$ of $SU(n)$, they establish a bijection
between non-trivial conjugacy classes in $G$  and so-called
crepant discrete valuations of the quotient variety
$\BC^n/G$.

Statements (3) and (4) were observed by Steinberg [{\bf St2}].
We note that there is an a priori proof by vonRandow [{\bf vR}]
that the group defined
by generators and relations in statement
 (8) of Theorem 4.1 is indeed
independent of the canonical ordering of the vertices of $\Delta$.
Equation (4-1) may be viewed as saying that the assignment
$v_i\mapsto g_i$ is a ``non-commutative harmonic map''.

\sec{5. Relation with the McKay correspondence}

Gonzalez-Sprinberg and Verdier [{\bf GS-V}] gave a geometric
construction of the McKay corerspondence in terms of the
first Chern class of some vector bundles on $\tilde X$.
Given any representation of $G$ in a a 
finite-dimensional vector space
$E$, there is a natural algebraic vector bundle
 $\calE$ over $\tilde X\min D
=X \min \{ o\}$. In terms of locally
 free sheaves, the locally free
sheaf $f_*\calO_{\BC^2\min \{ 0\}}$ has an action of $G$
hence it admits a decomposition into
 irreducible representations of $G$:

$$f_*\calO_{\BC^2\min \{ 0\}}=\oplus_{E\in Irr(G)}~E^{\pr} 
\otimes\calF_E\eqno(5-1)$$
where each $\calF_E$ is a locally free sheaf.
Then $\calE$ is the algebraic vector
 bundle corresponding to $\calF_E$.

We then want to extend the vector bundle $\calE$ to $\tilde X$.
For this it is enough to extend the locally free sheaf $\calF_E$
to $\tilde X$. Giving such an extension for any $E\in Irr(G)$
amounts to giving an extension of $f_*\calO_{\BC^2\min \{ 0\}}$
to a locally free sheaf over $\tilde X$. 
Let $j:\tilde X\min D\hook \tilde X$
be the inclusion. The extension given
in [{\bf GS-V}]  and [{\bf Kn}] is the 
sheaf $\calA$ of subalgebras of 
$j_*f_*\calO_{\BC^2\min \{ 0\}}$ generated
by $f_*\calO_{\BC^2}$ and by $\calO_{\tilde X}$.
It is proved in [{\bf GS-V}] [{\bf Kn}] that this
 is actually locally free.
The McKay correspondence is essentially 
given by the first Chern class
of the vector bundle $\calE$. We have the following

\prop{5.1}{(1) The group $H_2(\tilde X,\BZ)$ 
is the free abelian group
of rank $r$
generated by the homology classes $[D_i]$
 of the divisors $D_i$, for
$1\leq i\leq r$.
\item (2) The Picard group $Pic(\tilde X)$
 is isomorphic to 
$H^2(\tilde X,\BZ)$, hence to the $\BZ$-dual
 of $H_2(\tilde X,\BZ)=\BZ^r$.
The isomorphism $Pic(\tilde X)\to \BZ^r$ 
maps a line bundle $L$ to the
vector $(deg(L_{/D_i}))$.}

We identify both $H_2(\tilde X,\BZ)$ and 
$H^2(\tilde X,\BZ)$ with
$\BZ^r$. Let $(e_1,\cdots,e_r)$ denote 
the standard basis of
$H^2(\tilde X,\BZ)=\BZ^r$.

Let $\Lambda$
be the image of the natural  injection
Since $H_2(\tilde X,\BZ)$ identifies by Poincar\'e duality
with the cohomology group $H^2_c(\tilde X,\BZ)$ 
with compact supports,
there is a natural map $\kappa:H_2(\tilde X,\BZ)
\hook H^2(\tilde X,\BZ)$.

Then we have
\lem{5.2}{(1) The matrix
of $\kappa$ is the opposite of the Cartan matrix $C$.
\item (2) $\Lambda$ is a lattice in $H^2(\tilde X,\BZ)$.
Its index is equal to the connection index of the root
system.}

For each divisor $D_i$ the class $\kappa[D_i]
\in H^2(\tilde X,\BZ)$ will again be denoted by $[D_i]$.
It is the first Chern class of the locally free sheaf
$\calO(D_i)$.
It follows from Lemma 5.2 that the classes $[D_i]$ form a basis
of the $\BQ$-vector space $H^2(X,\BQ)$.
and that we have
$$e_i=-\sum_j~(C^{-1})_{ij}[D_j]\eqno(5-2)$$

The theorem of [{\bf GS-V}] can be stated as follows:

\theo{5.3}{(Gonzalez-Sprinberg [{\bf GS-V}], 
see also [{\bf Kn}])  
Let $E_i$ be the irreducible representation
of $G$ corresponding to the vertex $v_i$ of $\Delta$.
Then the first Chern class of the vector bundle $\calE_i$
is equal to the standard vector $e_i$ of $\BZ^r$.}

Theorem 5.3 indeed gives a geometric construction of the
McKay correspondence. The proofs  of the theorem  in
[{\bf GS-V}] and [{\bf Kn}] involve
a case by case computation, but there is a uniform proof in
[{\bf A-V}], which furthermore applies 
in arbitrary characteristic.

Now we will relate the vector bundles
$\calE$ to vector bundles with integrable connection.
We recall the notion of a Deligne vector bundle
with meromorphic integrable connection over an algebraic manifold
$Z$ with respect to divisor $Y\subset Z$ with normal
crossings. Let $V$ be a an algebraic vector bundle over $Z$.
Assume we have an integrable meromorphic 
connection $\na$ on $Z$, which is 
holomorphic over $Z\min Y$. Then we say that $(V,\na)$
is a {\it Deligne vector bundle with connection }if
\item(1) for any germ of holomorphic section $s$
of $V$, $\na(s)$ is a holomorphic section
of $\O^1_Z(log~Y)\otimes V$, where $\O^1_Z(log~Y)$
is the sheaf of $1$-forms with logarithmic poles
 (so $\na$ has at worst logarithmic poles
along $Y$);
\item(2) for any component $Y_j$ of $Y$, and any eigenvalue
$\a$ of the residue of $\na$ along $Y_j$, we have:
$$0\leq Re(\a)<1 \eqno(5-3)$$

The first Chern class of the vector bundle $V$ is computable
in terms of the residues of the connection along the components
$Y_j$ of $Y$. For each $j$, we have the cohomology class
$[Y_j]\in H^2(Z,\BC)$. The residue $Res_{Y_j}(\na)$ is a complex
number. Then we have:

\prop{5.4}{(Esnault-Verdier, see Appendix B of
 [{\bf E-V}]) Assume the vector bundle
$V$ with integrable meromorphic connection 
$\na$ satisfies (1) above.
Then we have:
$$c_1(V)=-\sum_j~Tr~Res_{Y_j}(\na)~[Y_j]
\in H^2(Z,\BC)\eqno(5-4)$$}

In fact, the result is proved in [{\bf E-V}] for a
 proper algebraic variety.
However, consider an algebraic vector bundle $V$ 
over $Z$ with integrable 
meromorphic connection satisfying (1). There
 exists a smooth compactification
$\bar X$ of $X$ such that $(\bar X\min X)
\cup \bar Y$ is a divisor with
 normal crossings. Then it follows from the
 theory of [{\bf De}] that $V$ can
 be extended to a vector bundle over $\bar X$ 
satisfying (1) with respect to
the divisor $(\bar X\min X)\cup \bar Y$. 
The equality (5-4) for $\bar V$ implies
the corresponding equality for $V$.

Here is an important class of examples of Deligne vector bundles
with integrable connection.

\lem{5.5}{Let $Y$ be a divisor with normal crossings in the smooth
complex algebraic variety $Z$. Let $S$ be a normal algebraic
variety and let $h:S\to Z$ be a proper morphism such that
\item (1) $h$ is an \'etale morphism over  $Z\min Y$.
\item (2) $h_*\calO_S$ is a locally free sheaf over $Z$.
\vskip .04 in
Then the vector bundle associated to $h_*\calO_S$ has the natural
structure of a Deligne line bundle with integrable connection.}

\proof{Over $Z\min Y$ we have a unique connection $\na$ 
on $h_*\calO_S$ which is compatible 
with the algebra structure.
For any section $u$ of $ h_*\calO_S$ over
 an open subset of $Z\min Y$,
 we can find a polynomial
equation 
$P(u)=0$, where $P(x)=x^n+\sum_{i=0}^{n-1}a_ix^i$ is a monic
 polynomial with coefficients in $\calO_Z$
such that $P^{\pr}(u)$ is nowhere vanishing.
Then we have:
$$\na(u)=\a\otimes {1\over P^{\pr}(u)},$$
where $\a=\sum_i~da_ix^i$. 
To prove properties (1)-(2) we work with holomorphic sheaves.
Now near a point of $Y$ where $Y$ has local equation
$x_1\cdots x_l=0$, the locally free sheaf $h_*\calO_S$
is a direct factor of the sheaf of algebras
$\calO_Z[x_1^{1\over m},\cdots,x_l^{1\over m}]$ for some $m$.
Indeed it is obtained as the subsheaf of invariants under some
subgroup of $(\mu_m)^l$. It is therefore enough to treat
the case of the sheaf of algebras
$\calO_Z[x_1^{1\over m},\cdots,x_l^{1\over m}]$.
This has a basis consisting of functions of the type
$u=x_1^{q_1\over m}\cdots x_l^{q_l\over m}$
for $0\leq q_j\leq m-1$. Then we have
$$\na u=\sum_j~{q_j\over m}{dx_j\over x_j}\otimes u,$$
which makes (1) and (2) apparent. \qed}

This however does not directly apply to our
 bundle of algebras
 over $\tilde X$, because the corresponding 
 sheaf of algebras $\calA$
is not integrally closed. We need to introduce the integral
closure $\tilde\calA$. This satisfies all the assumptions
of proposition 5.3, at least on $\tilde X\min T$, where
$T\subset \tilde X$ is a finite set. Since we are interested
in the first Chern class, we may neglect the effect of deleting
this finite set. We then have an exact sequence of $G$-equivariant
coherent
sheaves on $\tilde X$:
$$0\to \calA\to\tilde\calA\to \calF\to 0\eqno(5-5)$$
for some sheaf $\calF$ supported on $Y$.
For the isotypic components associated to $E\in Irr(G)$ this
yields an exact sequence
$$0\to \calE \to\tilde\calE\to \calH\to 0\eqno(5-6)$$
for some coherent sheaf $\calH$ supported on $Y$.
For the first Chern classes we obtain
$$c_1(\calE)=c_1(\tilde\calE)-\sum_{j=1}^rm_j[D_j]\eqno(5-7)$$
where $m_j\geq 0$.

We can state

\prop{5.6}{Let $\calE$ be the vector bundle
over $\tilde X$ associated to an irreducible
representation $E$ of $G$. We have
$$c_1(\calE)=-\sum_{j=1}^r~q_j[D_j]~{\rm with~}q_j
\geq 0\eqno(5-8)$$}

\proof{The vector bundle $\tilde E$ has an integrable
connection, and is a direct factor of the vector bundle
with integrable connection $\calA$. We have from Proposition 5.4 the equality
$$c_1(\tilde E)=-\sum_j \lambda_j [D_j],$$
where $\lambda_j$ is the trace of the residue of the integrable
connection along $D_j$. By Lemma 5.5 all the eigenvalues
of the residue have real part $\geq 0$, so that $Re(\lambda_j)\geq 0$.
Thus we get the equality
$$c_1(\tilde E)=-\sum_j (\lambda_j+q_j) [D_j],$$
with $Re(\lambda_j+q_j)\geq 0$.
Since the $[D_j]$ are linearly independent, it follows
that the $\lambda_j+q_j$ are positive rational numbers. \qed}

On the other hand we know from Theorem 5.3 that the $q_j$ are
(up to sign) equal  to the
coefficients of the inverse of the Cartan matrix. Therefore
we conclude

\cor{5.7}{The coefficients $(C^{-1})_{ij}$
of the inverse $C^{-1}$ of the Cartan matrix are all
$\geq 0$.}

This can of course be easily read off the tables
in Bourbaki [{\bf Bo}]. Indeed $(C^{-1})_{jk}$
is the coefficient of $\a_j$ in the fundamental weight
$\omega_k$. This was pointed out to me by Dolgachev.
We note the graph-theoretic
interpretation of $C^{-1}$ given  by Lusztig and Tits [{\bf L-T}].

It is actually easy to prove directly
 that all coefficients of 
$C^{-1}$ are positive. Recall that $C=2I-A$ 
where $A$ is the matrix
of $\Delta$. Now the operator norm of $A$ 
is the norm of the largest
eigenvalue, which is well-known to be
 smaller than $2$. Thus we have 
a convergent series:
$$C^{-1}=(2I-A)^{-1}={1\over 2}
\sum_{n\geq 0}~2^{-n}(A^n)_{ij}\eqno(5-9)$$
The coefficient $(A^n)_{ij}$ is the number of paths of
 length $n$ from $v_i$ to $v_j$
in $\Delta$ . So it is always $\geq 0$ and given
$(i,j)$ there exists $n$ such that $(A^n)_{ij}>0$.
In fact we have an estimate:
$$(C^{-1})_{ij}\geq {2^{1-n}\over 3}\eqno(5-10)$$
where $n$ is the distance between the vertices $v_i$ and $v_j$.
This estimate is sharp for the $A_2$ case.

Our final result involves both the  McKay correspondence
and the dual correspondence.
Let $g_j$ be a representative of the conjugacy class
associated to the divisor $D_j$. Then we have:

\prop{5.8}{For the irreducible representation $E_k$ of $G$
associated to the vertex $v_k$ of $\Delta$
we have
$$det(g_j,E_k)=exp(-2\pi i(C^{-1})_{jk})\eqno(5-11)$$}

\proof{The conjugacy class of the monodromy operator
$g_j$ on $E_k$ is represented by the operator
$exp(-2\pi iRes_{D_j}\na)$, where $\na$ is the
meromorphic connection on the Deligne vector bundle
$\tilde\calE_k$; this is a general property of the residue
in the case of semisimple monodromy,
proved in [{\bf De}]. So we have:
$$det(g_j,E_k)=exp(-2\pi i Tr~Res_{D_j}\na).$$
 Now by Proposition 5.4  the trace of the residue
is the opposite of the coefficient of $c_1(L_j)$
in $c_1(\tilde\calE_k)$. This coefficient is congruent
modulo $\BZ$ to the similar coefficient
for $c_1(\calE_k)$. The latter coefficient
is the opposite of the coefficient $(C^{-1})_{jk}$ of $C^{-1}$. \qed}

Recall that the connection index of $\Delta$ is equal to the determinant
of the Cartan matrix.

\cor{5.9}{The exponent of the abelianization $G^{ab}$ of $G$
is equal to the connection index of $\Delta$.}

\proof{Indeed if $n$ is the connection index, then all coefficients
$(C^{-1})_{ij}$ belong to ${1\over n}\BZ$. It then follows from
(5-11) that any character $G\to\BC^*$ takes 
values in the $n$-th roots of unity.
So the exponent $m$ of $G^{ab}$ divides $n$. Now let $n=mq$; i
it follows from (5-11) that for each coefficient  $(C^{-1})_{ij}$
we have $m(C^{-1})_{ij}\in\BZ$. But it is easy to see that the g.c.d
of the integers $n(C^{-1})_{ij}$ is equal to $1$. This implies $q=1$
and $m=n$. \qed}

This result in fact follows easily from the presentation of the group
$G$ in terms of the Cartan matrix given in [{\bf H-N-K}].

For instance, the connection index of the diagram $E_8$ is equal to
$1$, which implies that the binary icosahedral group is perfect
(a well-known fact, of course).

This gives evidence for the idea that the matrix-valued
Fourier transform obtained by combining the two types
of correspondences should be very significant geometrically.
One easily checks in the $A_n$ case that $(C^{-1})_{jk}$
is congruent to $-jk\over {n+1}$ modulo $\BZ$.
Therefore we get an automorphism $F$
of the space of functions on $G$, whose matrix
is 
$$F_{jk}=exp({-2\pi i jk\over n+1})\eqno(5-12).$$
This is of course the usual Fourier transform on the cyclic
group $\mu_{n+1}$.
For $G$ non-abelian, we obtain a matrix-valued Fourier
transform; by taking the trace of a representation, we
obtain an automorphism of the space of central functions
on $G$; however, this is not of finite order, already
in the case of the binary octahedral group (case $D_4$). 

This strongly suggests that the main object of interest should
be the matrix-valued Fourier transform, not just 
the scalar-valued Fourier transform.

\vskip 1pc

REFERENCES
\vskip .5pc

[{\bf A}] M.Artin, \it On isolated rational singularities
of surfaces, \rm Amer. J. Math. \bf 88 \rm (1966), 129-136

[{\bf A-V}] M. Artin and J-L. Verdier, \it Reflexive modules
over rational double points, \rm Math. Ann. \bf 270 \rm (1985), 79-82

[{\bf Bo}] N. Bourbaki, \it Lie Groups and Lie Algebras,
Chapters 4,5,6

[{\bf Br}] E. Brieskorn, \it
Rationale singularit\"aten komplexer Fl\"achen,
\rm Invent. Math. \bf 4 \rm(1968), 336-358

[{\bf Co}] H. S. M. Coxeter, \it
Regular Complex Polytopes, \rm Cambridge Univ. Press
(1974)

[{\bf De}] P. Deligne,  \it Equations Diff\'erentielles
\`a Points Singuliers R\'eguliers, \rm
Lecture Notes on Math vol. \bf 163 \rm (1970), \Sp 

[{\bf Du Val1}] P. Du Val, \it
On isolated singularities which do not affect the condition
of adjunction, \rm Proc. Cambridge Phil. Soc. \bf 30
\rm (1934), 453-465

[{\bf Du Val2}] P. Du Val, \it
Homographies, Quaternions and Rotations, \rm Clarendon Press
(1964)

[{\bf E-V}] H. Esnault and E. Viewhweg, \it
Logarithmic de Rham complexes and vanishing theorems, \rm
Invent. Math. \bf 86 \rm (1986), 161-194

[{\bf GS-V}] G. Gonzalez-Sprinberg and J-L. Verdier,
\it Construction g\'eom\'etrique de la correspondence
de McKay, \rm Ann. Sc. ec. Norm. Sup. \bf 16
\rm (1983), 409-449

[{\bf H1}] F. Hirzebruch, \it
\"Uber vierdimensionale Riemmansche Fl\"achen
mehr-\hfill\break
-deutiger analytischer Funktionen von zwei
komplexer Ver\"anderlichen, \rm Math. Ann. \bf 126
\rm (1953), 1-22

[{\bf H-N-K}]  F. Hirzebruch, W. D. Neumann
and S. S. Koch, \it
Differentiable manifolds and Quadratic Forms, \rm
Lecture Notes in Pure and Applied Math. vol \bf 4 \rm,
Marcel Dekker (1971)

[{\bf I-R}] Y. Ito and M. Reid, \it The Mckay
correspondence for finite subgroups of $SL(3,\BC)$,\rm
\hfill\break
preprint (1994), al-geom/9411010

[{\bf Kn}] H. Kn\"orrer, \it Group representations and
the resolution of rational double points, \rm Contemp. Math. 
vol. \bf 45 \rm 91985), 175-221

[{\bf Ko}] B. Kostant, \it On finite subgroups of
$SU(2)$, simple Lie algebras and the McKay correspondence,
\rm Ast\'erisque vol. Hors-S\'erie (1985), 109-255

[{\bf L-T}] G. Lusztig and J. Tits, \it The inverse of a Cartan
matrix, \rm Ann. Univ. Timi\c soara Ser \c Stiin\c t. Mat. \bf 30
\rm (1992), no. 1, 17-23

[{\bf McKay1}] J. McKay, \it Graphs, singularities and
finite groups, \rm Proc. Symp. Pure Math. \bf 37 \rm (1980),
183-186 

[{\bf McKay2}] J. McKay, \it Cartan matrices, finite groups
of quaternions, and \hfill\break
kleinian singularities, \rm
Proc. Amer. Math. Soc. \bf \rm (1981), 153-154

[{\bf Mu}] D. Mumford, \it
The topology of  normal singularities of an algebraic surface
and a criterion for simplicity, \rm Publ. Math. IHES \bf 9 \rm 
(1961), 23-64

[{\bf St1}] R. Steinberg, \it Kleinian singularities
and unipotent elements, \rm
Proc. Symp. Pure Math. \bf 37 \rm (1980), 265-270

[{\bf St2}] R. Steinberg, \it
Subgroups of $SU_2$, Dynkin diagrams and affine Coxeter elements,
\rm Pacific Jour. Math. \bf 118 \rm (1985), 587-598

[{\bf vR}] R. von Randow, \it Zur Topologie von dreidimensionalen
Baummannifal-\hfill\break
-tigkeiten, \rm Bonner Math. Schriften \bf 14
\rm (1962)
\vskip 1.5pc
Penn State University

Department of Mathematics

305 McAllister

University Park, PA. 16802

USA

\vskip .5pc
e-mail:jlb@math.psu.edu

\bye